\documentclass[12pt]{article}
\usepackage{amssymb,amsmath,amsthm}
\usepackage[dvips]{graphicx,xcolor}

\newcommand{\pf}{\vspace{3mm}\noindent{\it Proof)}\hspace{5mm}}

\newcommand{\R}{{\bf R}}

\newcommand{\p}{\partial}

\newcommand{\ep}{\varepsilon}

\newcommand{\al}{\alpha}
\newcommand{\bt}{\beta}

\newcommand{\dsp}{\displaystyle}
\newcommand{\lr}{\longrightarrow }

\newcommand{\fa}{\frac}

\newcommand{\lp}[1]{\left(#1 \right)}
\newcommand{\lb}[1]{\left\{#1 \right\}}

\begin{document}
 
\addtocounter{footnote}{1}
 \title{}
 
\begin{center}
{\large\bf On heatlike lifespan of solutions of semilinear wave equations \\
in Friedmann-Lema\^itre-Robertson-Walker spacetime\\
 
}
\end{center}
\vspace{3mm}
\begin{center}

Kimitoshi Tsutaya$^\dagger$ 
        and Yuta Wakasugi$^\ddagger$ \\
\vspace{1cm}

 $^\dagger$Graduate School  of Science and Technology \\
Hirosaki University  \\
Hirosaki 036-8561, Japan\\
\footnotetext{* The research was supported by JSPS KAKENHI Grant Number JP18K03351. }

\vspace{5mm}
        $^\ddagger$
Graduate School of Engineering \\
Hiroshima University \\
Higashi-Hiroshima, 739-8527, Japan

\vspace{5mm}
{\it In memory of Prof. K\^oji Kubota}

\end{center}

\vspace{3mm}
\begin{abstract}
Consider a nonlinear wave equation  for a massless scalar field  with self-interaction 
in the spatially flat Friedmann-Lema\^itre-Robertson-Walker spacetimes. 
We treat the so-called heatlike case where the critical exponent is affected by the Fujita exponent. 
We show upper bounds of the lifespan of blow-up solutions by distinguishing subcritical and critical cases. 

\end{abstract}  



\section{Introduction.}
\addtolength{\baselineskip}{2.4mm}
\quad 
This paper is subsequent to our recent work \cite{TW1} concerned with  
the semilinear wave equation in the Friedmann-Lema\^itre-Robertson-Walker (FLRW) spacetimes. 
The spatially flat FLRW metric is given by 
\[
g: \; ds^2=-dt^2+a(t)^2d\sigma^2, 
\]
where the speed of light is equal to $1$, $d\sigma^2$ is the line element of $n$-dimensional Euclidean space and $a(t)$ is the scale factor, which describes expansion or contraction of the spatial metric. 
If we solve the Einstein equation with the energy-momentum tensor for the perfect fluid under an equation of state, we obtain 
\begin{equation}
a(t)=ct^{\fa 2{n(1+w)}}  
\label{scw}
\end{equation}
with some constant $c$, where $w$ is a proportional constant in the equation of state. See \cite{CGLY} for details.  The present paper treats the case $2/n-1<w \le 1$ for $n\ge 2$, that is, a decelerating expanding universe. 

For the flat FLRW metric with \eqref{scw},  
the semilinear wave equation $\Box_g u\\
=|g|^{-1/2}\p_\al(|g|^{1/2}g^{\al\bt}\p_\bt)u= -|u|^p$ with $p>1$ becomes 
\begin{equation}
u_{tt}-\fa 1{t^{4/(n(1+w))}}\Delta u+\fa 2{(1+w)t}u_t=|u|^p,   \quad  x\in \R^n 
\label{ore}
\end{equation}
where  $\Delta=\p_1^2+\cdots \p_n^2, \; \p_j=\p/\p x^j, \; j=1,\cdots,n, \; (x^1,\cdots, x^n)\in \R^n$
.  
Our aim of this paper is to show that blow-up in a finite time occurs for the equation above as well as upper bounds of the lifespan of the blow-up solutions. 

In order to compare with the related known results including the case of the Minkowski spacetime, we first consider the Cauchy problem
\begin{equation}
\begin{cases}
\dsp
u_{tt}-\fa 1{t^{2\al}}\Delta u+\fa\mu{t}u_t=|u|^p, &\qquad  t>1, \; x\in \R^n  \\
  & \\
u(1,x)=\ep u_0(x), \; u_t(1,x)=\ep u_1(x), &\qquad x\in \R^n,  
\end{cases}
\label{Prob0}
\end{equation}
where $\al$ and $\mu$ are nonnegative constants, and $\ep>0$ is a small parameter.  
We then return to \eqref{ore}.  

Let $T_\ep$ be the lifespan of solutions of \eqref{Prob0}, that is, $T_\ep$ is the supremum of $T$ such that  \eqref{Prob0} have a solution for $x\in \R^n$ and $1\le t<T$. 
Let $p_F(n)$ denote the Fujita exponent $p_F(n)=1+2/n$, and $p_S(n)$ the Strauss exponent which is the positive root of the equation 
\begin{equation}
\gamma_S(n,p)\equiv -(n-1)p^2+(n+1)p+2=0. 
\label{snp}
\end{equation}
For the special case $\al=0$, there have been many results. 
Among them, it is shown by Tu and Lin \cite{LinTu,TuLin1} and Ikeda and Sobajima \cite{IkSo18} that 
\begin{align}
&T_\ep \le C\ep^{-2p(p-1)/\gamma_S(n+\mu,p)}\quad \mbox{ if }\quad  1< p < p_s(n+\mu), \; n\ge 2, 
\label{subst}\\
&T_\ep \le \exp(C\ep^{-p(p-1)})\quad  \mbox{ if }\quad p = p_s(n+\mu)>p_F(n), \; n\ge 2.  \nonumber
\end{align}
For the case $\al=2/3, \; \mu=2$ and $n\ge 1$ in \eqref{Prob0}, 
it is proved by Galstian and Yagdjian \cite{GY} that finite time blow-up occurs if $1<p<p_F(n/3)$ or $1<p<p_{cr}(n)$, where $p_{cr}(n)$ is the positive root of 
\[
-(n+3)p^2+(n+13)p+2=0. 
\] 
Finite time blow-up for $0<\al<1$ and $\mu=2$ is also shown in \cite{GY}.  

The authors have recently proved blow-up in a finite time and upper estimates of the lifespan for the case $0\le \al <1$. 
Let 
\begin{align}
\gamma(n,p,\al,\mu)=-p^2\lp{n-1+\fa{\mu-\al}{1-\al}}+p\lp{n+1+\fa{\mu+3\al}{1-\al}}+2
\label{gm} 
\end{align}
and $p_c(n,\al,\mu)$ be the positive root of the equation $\gamma(n,p,\al,\mu)=0$. 
It is shown in \cite{TW1} that 
\begin{align}
&T_\ep\le C\ep^{\fa{-2p(p-1)}{(1-\al)\gamma(n,p,\al,\mu)}} && \mbox{ if }\quad 1<p<p_c(n,\al,\mu), \label{Stsub} \\
&T_\ep\le C\ep^{-\fa{p-1}{2-\{n(1-\al)+\mu-1\}(p-1)}} &&\mbox{ if }\quad p<1+2/(n(1-\al)+\mu-1), \label{upb0}  \\
&T_\ep\le \exp(C\ep^{-p(p-1)}) && \mbox{ if }\quad p=p_c(n,\al,\mu)>p_F(n(1-\al))=1+2/(n(1-\al)). \nonumber
\end{align}
Recently, Palmieri \cite{Pa} has independently obtained similar results to above 
by using different methods. 

We remark that if $\mu=\al=0$, the $p_c(n,0,0)$ coincides with the Strauss exponent $p_S(n)$. 
If $\al=2/3$ and $\mu=2$, our blow-up range of $p$ is the same as $1<p<p_{cr}(n)$ in \cite{GY}. 
We note that these upper bounds stated so far are sharp 
if the power $p$ is dominated by the Strauss exponent or $p_c(n,\al,\mu)$. 
We can say in this case that the lifespan of solutions is {\it wavelike} as mentioned in \cite{LTW}.  

On the other hand, if the critical exponent for \eqref{Prob0} is affected by the Fujita exponent, the upper bound of the lifespan is expected to become sharper. 
In fact if $\al=0$, $p< p_F(n)$ and $p<2/(n-\mu+1)$,  
the following upper bound of the lifespan is sharper than \eqref{subst}:  
\[
T_\ep \le
 C\ep^{-(p-1)/\{2-n(p-1)\}}. 
\]
The estimate above holds for $1<p< p_F(n)$ and $\mu>1$. 
We can say that this case 
 corresponds to the {\it heatlike} lifespan as in \cite{LTW}. 
For the critical case, there holds
\[
T_\ep \le
\exp(C\ep^{-(p-1)})  \quad \mbox{ if \quad $1<p= p_F(n)$ and $\mu>1$}. 
\]
These upper bounds are proved by the second author \cite{Wa1,Wa2} and Ikeda et al \cite{ISW}.  

In this paper we treat the case where the critical exponent for \eqref{Prob0} with $0\le \al <1$ is affected by $p_F$, that is, the heatlike case, and show upper estimates of the lifespan of the blow-up solutions. 
We then apply our results to the original equation \eqref{ore}. 
Our aim is especially to clarify the difference with the case of the Minkowski spacetime and also how the scale factor affects the lifespan of the solution.  

The following theorem is our main result.

\noindent
{\bf Theorem 1.1.} 
{\it 
Let $n\ge 2, \; 0\le \al<1, \; \mu\ge 0$ and $1<p\le p_F(n(1-\al))$.
Assume that $u_0\in C^2(\R^n)$ and $u_1\in C^1(\R^n)$ are nonnegative, nontrivial and $\mbox{\rm supp }u_0, \mbox{\rm supp }u_1\subset \{|x|\le R\}$ with $R>0$. 
Suppose that the problem \eqref{Prob0} has  a classical solution $u\in C^2([1,T)\times\R^n)$.  
Then $T<\infty$ and there exists a constant $\ep_0>0$ depending on 
$p,\al,\mu,R,u_0,u_1$ such that $T_\ep$ has to satisfy 
\begin{align*}
T_\ep&\le C\ep^{-\fa{p-1}{2-n(1-\al)(p-1)}} && (p<p_F(n(1-\al))) \\
T_\ep&\le \exp\lp{C\ep^{-p(p-1)/(p+1)}} &&(p=p_F(n(1-\al)) \mbox{ and }0\le \mu \le 1) 
\\
T_\ep&\le \exp\lp{C\ep^{-(p-1)}} &&(p=p_F(n(1-\al)) \mbox{ and }\mu > 1) 
\end{align*}
for $0<\ep\le \ep_0$, 
where $C>0$ is a constant independent of $\ep$.
}

\noindent
{\bf Remark.}\quad (i) In the critical case $p=p_F(n(1-\al))$, if $0\le \mu \le 1$, then the estimate above is better than that for the case $\mu >1$. 
However, \eqref{Stsub} and \eqref{upb0} are applicable for $p=p_F(n(1-\al))$ and 
$0\le \mu \le 1$. These estimates are much better than $T_\ep\le \exp\lp{C\ep^{-p(p-1)/(p+1)}}$ above.   
See Figure 1 in Section 4. \\
(ii) Similar results are independently shown by \cite{Pa} where energy solutions are treated.

\vspace{5mm}    
The paper is organized as follows. 
In Section 2 we show the theorem in the subcritical case on the power $p$ by the generalized Kato's lemma proved in \cite{TW1}. 
In Section 3 we provide another Kato's lemma as a new tool to show the theorem in the critical case. 
In Section 4 we apply the theorem to the original equation \eqref{ore}, and discuss the effect of the scale factor to the solutions.

\vspace{5mm}

\section{Subcritical case.} 
\setcounter{equation}{0}

\quad We prove Theorem 1.1 in the subcritical case $p<p_F(n(1-\al))$.
We use the following lemma, which is a generalized Kato's lemma, to prove Theorem 1.1.

\noindent
{\bf Lemma 2.1}
{\it Let $p>1, \;a\ge 0, \;b>0, \;q>0, \; \mu\ge 0$ and 
\[
M\equiv (p-1)(b-a)-q+2>0. 
\]
Let $T\ge T_1>T_0\ge 1$. 
Assume that $F\in C^2([T_0,T))$ satisfies the following three conditions:  
\begin{align*}
(i) \quad & F(t) \ge A_0t^{-a}(t-T_1)^b \qquad \mbox{for }t\ge T_1 \\
(ii) \quad &F''(t) +\fa{\mu F'(t)}{t}\ge A_1(t+R)^{-q}|F(t)|^p \quad \mbox{for }t\ge T_0	 \\
(iii) \quad &F(T_0)\ge 0, \quad F'(T_0)>0,  
\end{align*}
where $A_0, A_1$ and $R$ are positive constants. Then 
$T$ has to satisfy
\[
T<C A_0^{-(p-1)/M},  
\]
where $C$ is a constant depending on $R,A_1,\mu,p,q,a$ and $b$. 
}

See \cite{TW1} for its proof.

\vspace{1cm}

We now prove Theorem 1.1 by applying Lemma 2.1. 
Set 
\[
F(t)=\int u(t,x)dx.
\] 
By finite speed of propagation, we have 
\begin{equation}
\mbox{supp }u(t,\cdot)\subset \{|x|\le A(t)+R\}, 
\label{suppu}
\end{equation}
where 
\[
A(t)=\int_1^t s^{-\al}ds =\fa{t^{1-\al}-1}{1-\al}. 
\]
See Appendix in \cite{TW1} for its proof. 
Integrating the equation \eqref{Prob0} and using H\"older's inequality, we have by \eqref{suppu}, 
\begin{align}
F''(t)+\fa\mu t F'(t) &=\int |u|^p dx  
   \label{eqF} \\
 & \ge \fa 1{(A(t)+R)^{n(p-1)}}|F(t)|^p \nonumber \\
 & \ge \fa C{(t+R)^{n(1-\al)(p-1)}}|F(t)|^p. 
\label{Ho}
\end{align}
Mutiplying \eqref{Ho} by $t^\mu$
 and integrating imply
\begin{equation}
t^\mu F'(t)-F'(1)\gtrsim \int_1^t s^{\mu-n(1-\al)(p-1)} |F(s)|^pds\ge 0. 
\label{F1}
\end{equation}
Note that 
\[
F'(t)>0 \quad \mbox{for }t\ge 1
\]
since $F'(1)>0$ by assumption. 
Moreover, from \eqref{F1}, 
\begin{align*}
F'(t)&\gtrsim t^{-\mu}F'(1). \\
\intertext{Integrating over $[1,t]$ implies}
F(t)-F(1)&\gtrsim F'(1)\int_1^t s^{-\mu}ds >0.
\end{align*}
We hence see that 
\begin{equation} 
F(t)\ge F(1)=C\ep>0 \qquad \mbox{for }t\ge 1
\label{F2}
\end{equation}
by assumption. 
From \eqref{Ho} and \eqref{F2}, 
we have 
\begin{equation}
\int |u|^pdx \gtrsim \ep^p t^{-n(1-\al)(p-1)} \qquad \mbox{for }t\ge 1. 
\label{up}
\end{equation}
On the other hand, mutiplying \eqref{eqF} by $t^\mu$
 and integrating imply
\begin{align*}
t^\mu F'(t)-F'(1)&= \int_1^t s^\mu\int |u|^p dxds. 
\intertext{Since $F'(1)>0$ by assumption,}  
 F'(t)& \ge t^{-\mu}\int_1^t s^\mu\int |u|^p dxds. 
\end{align*} 
Integrating again, we have from $F(1)>0$ by assumption, 
\begin{equation}
F(t)\ge \int_1^t \tau^{-\mu}\int_1^\tau s^\mu\int |u|^p dxdsd\tau  \qquad 
\mbox{for }t\ge 1. 
\label{Fup}
\end{equation}  
Using \eqref{Fup} implies 
\begin{align}
F(t)&\ge C\ep^p \int_1^t \tau^{-\mu}\int_1^\tau s^{\mu-n(1-\al)(p-1)}dsd\tau \nonumber \\
&\ge C\ep^p \int_1^t \tau^{-\mu-n(1-\al)(p-1)}
  \int_1^\tau (s-1)^\mu dsd\tau \nonumber\\
&\ge C\ep^p  \int_1^t \tau^{-\mu-n(1-\al)(p-1)}(\tau-1)^{\mu+1}d\tau \nonumber\\ 
&\ge C\ep^p t^{-\mu-n(1-\al)(p-1)}\int_1^t (\tau-1)^{\mu+1}d\tau.  \nonumber \\
\intertext{Therefore, we obtain}
F(t)&\ge C\ep^p t^{-\mu-n(1-\al)(p-1)}(t-1)^{\mu+2} 
\qquad \mbox{for }t\ge 1. 
\label{ann}
\end{align}

Finally, by \eqref{Ho} and \eqref{ann}, applying Lemma 2.1 with $q=n(1-\al)(p-1)$, $a=\mu+n(1-\al)(p-1)$, $b=\mu+2$ and $A_0=C\ep^p$, we obtain the desired results since
\begin{align*}
M &=(p-1)\lb{2-n(1-\al)(p-1)}-n(1-\al)(p-1)+2\\
&=p\lb{2-n(1-\al)(p-1)}>0. 
\end{align*}
This completes the proof of Theorem 1.1 in the subcritical case. 
\hfill\qed \\

\section{Critical case.} 
\setcounter{equation}{0}
\quad We prove the theorem for the critical case $p=p_F(n(1-\al))$, i.e., $n(1-\al)(p-1)-2=0$ in this section. 
We first show the following lemma, which is another generalized Kato's lemma.

\noindent
{\bf Lemma 3.1}
{\it \quad 
Let $p>1, \;b>0, \; \mu\ge 0$ and $T\ge T_1>T_0\ge 1$. 
Assume that $F\in C^2([T_0,T))$ satisfies the following three conditions:  
\begin{align*}
(i) \quad & F(t) \ge A_0\lp{\ln \fa t{T_1}}^b\qquad \mbox{for }t\ge T_1, \\
(ii) \quad &F''(t) +\fa{\mu F'(t)}{t}\ge A_1(t+R)^{-2}|F(t)|^p \quad \mbox{for }t\ge T_0,	 \\
(iii) \quad &F(T_0)\ge 0, \quad F'(T_0)>0,  
\end{align*}
where $A_0, A_1$ and $R$ are positive constants. Then,  
$T$ has to satisfy
\[
T<
\begin{cases}
\exp\lp{CA_0^{-(p-1)/\{b(p-1)+2\}}} & \mbox{if }\mu\le 1, \\
\exp\lp{CA_0^{-(p-1)/\{b(p-1)+1\}}} & \mbox{if }\mu > 1, 
\end{cases}
\] 
where $C$ is a constant depending on $R,A_1,\mu,p$ and $b$. 
}

\pf 
We first integrate assumption (ii) multiplied by $t^\mu$ over $[T_0,t]$. By assumption (iii) and integration again, we have 
\begin{align}
F(t)& \ge A_1C_R \int_{T_0}^t \tau^{-\mu}\int_{T_0}^\tau s^{\mu-2}|F(s)|^pdsd\tau.  \label{1}
\end{align}
See \cite{TW1}. 
We divide the proof into two cases. 
\\
(i) $\mu\le 1$. From \eqref{1} and assumption (i), we have
\begin{align*}
F(t) &\ge A_0^pA_1C_R \int_{T_1}^t \tau^{-\mu}\int_{T_1}^\tau s^{\mu-2}\lp{\ln\fa s{T_1}}^{bp}dsd\tau  \\
 &\ge 
\fa{A_0^pA_1C_R}{bp+1}\int_{T_1}^t \tau^{-1}\lp{\ln\fa \tau{T_1}}^{bp+1}d\tau \\
&= \fa{A_0^pA_1C_R}{(bp+1)(bp+2)}\lp{\ln\fa t{T_1}}^{bp+2}
\qquad \mbox{for }t\ge T_1. 
\end{align*} 
Based on the fact above, we define the sequences $b_j$ and $C_j$ for $j = 0, 1, 2, \cdots$ by
\begin{eqnarray}
b_{j+1} = pb_j+2, \quad & C_{j+1} = \dsp\frac{A_1C_RC_j^p}{(pb_j+2)^2}  \label{e2}\\
b_0 = b, \quad & C_0 = A_0.  \label{e3}
\end{eqnarray}
Solving \eqref{e2} and \eqref{e3}, we obtain 
\[
b_j = p^j\lp{b+\fa{2}{p-1}}-\fa{2}{p-1},
\]
and thus,
\[
C_{j+1} = \frac{A_1C_RC_j^p}{b_{j+1}^2} \ge \lp{b+\fa{2}{p-1}}^{-2}\frac{A_1C_RC_j^p}{p^{2j+2}}.
\]
Then, 
\begin{align*}
C_j &\ge \frac{BC_{j-1}^p}{p^{2j}}\\
&\ge \fa B{p^{2j}}\lp{\fa{BC_{j-2}^p}{p^{2(j-1)}}}^p =\fa{B^{1+p}}{p^{2j+2p(j-1)}}C_{j-2}^{p^2} \\
&\ge \cdots\cdots \ge \fa{B^{1+p+p^2+\cdots+p^{j-1}}}{p^{2(j+p(j-1)+p^2(j-2)+\cdots+p^{j-1})}}C_0^{p^j}, \\
\intertext{and }
\ln C_j &\ge \ln B \sum_{k=0}^{j-1}p^k-2\ln p \sum_{k=0}^j kp^{j-k}+p^j\ln C_0 
\\
&= \frac{\ln B}{p-1}(p^j-1)-2p^j \sum_{k=0}^j 
\frac{k}{p^k}\ln p+p^j\ln C_0, 
\end{align*}
where $B=\lp{b+\fa{2}{p-1}}^{-2}A_1C_R$. 
For sufficiently large $j$, we have
\[
C_j \ge \exp(Ep^j),
\] 
where
\[
E = \frac{1}{p-1}\min\left(0, \ln B \right) - 2\sum_{k=0}^{\infty}\frac{k}{p^k}\ln p+\ln A_0. 
\]
Thus, since $F(t)\ge C_j\lp{\ln\fa t{T_1}}^{b_j}
$ holds for $t\ge T_1$, we obtain
\begin{equation}
F(t) \ge \lp{\ln\fa t{T_1}}^{-2/(p-1)} 
 \exp\left[\left\{E+\lp{b+\fa{2}{p-1}}\ln\lp{\ln\fa t{T_1}}\right\}p^j\right]
  \qquad \mbox{for }t\ge T_1. 
\label{5}
\end{equation}
By choosing $t$ large enough, we can find a positive $\delta$ such that
\[
E+\lp{b+\fa{2}{p-1}}\ln\lp{\ln\fa t{T_1}}\ge \delta > 0. 
\]
It then follows from \eqref{5}  that $F(t) \lr \infty$ as $j \to \infty$ for 
sufficiently large $t$. We therefore see that the lifespan $T$ of $F(t)$ has to satisfy
\[
T<\exp\lp{CA_0^{-(p-1)/\{b(p-1)+2\}}},  
\]
where $C$ is a constant depending on $R,A_1,\mu,p$ and $b$. \\
(ii) $\mu >1$. 
We first introduce a bounded and increasing sequence following \cite{AKT}. 
Let $a_j=1+\fa 12+\cdots +\lp{\fa 12}^j$. Note that 
\begin{align}
&\fa{a_j}{a_{j+1}}=\fa{a_j}{1+\fa 12a_j}\ge \fa{a_j}{a_j+\fa 12a_j}=\fa 23, \label{6}\quad \intertext{and }\quad 
&1-\fa{a_j}{a_{j+1}}=\fa{\lp{\fa 12}^{j+1}}{1+\fa 12+\cdots +\lp{\fa 12}^j}>
\lp{\fa 12}^{j+2} \quad \mbox{for }j=0,1,2,\cdots. \label{7}
\end{align}
To avoid repeating the same argument, we start with the assumption
\[
F(t) \ge A_0\lp{\ln \fa t{a_jT_1}}^b \qquad \mbox{for }t\ge a_jT_1 \mbox{ and }
j=0,1,2,\cdots.
\]
If $j=0$, this is the same as assumption (i). 
From \eqref{1}, 
\[
F(t)\ge A_0^pA_1C_R \int_{a_jT_1}^t \tau^{-\mu}\int_{a_jT_1}^\tau s^{\mu-2}
\lp{\ln\fa s{a_jT_1}}^{bp}dsd\tau. 
\]
For $t\ge a_{j+1}T_1$, we have 
\begin{align*}
F(t)&\ge A_0^pA_1C_R \int_{a_{j+1}T_1}^t \tau^{-\mu}\int_{a_j\tau/a_{j+1}}^\tau s^{\mu-2}\lp{\ln\fa s{a_jT_1}}^{bp}dsd\tau \\
&\ge A_0^pA_1C_R\lp{\fa{a_j}{a_{j+1}}}^\mu\lp{1-\fa{a_j}{a_{j+1}}} 
\int_{a_{j+1}T_1}^t \tau^{-1} \lp{\ln\fa \tau{a_{j+1}T_1}}^{bp}d\tau. 
\end{align*}
Using \eqref{6} and \eqref{7} yields 
\begin{align*}
F(t)&\ge \fa{A_0^pA_1C_R}{2^{j+2}}\lp{\fa 23}^\mu 
\int_{a_{j+1}T_1}^t \tau^{-1} \lp{\ln\fa \tau{a_{j+1}T_1}}^{bp}d\tau \\
&=\fa{A_0^pA_1C_R}{2(bp+1)2^{j+1}}\lp{\fa 23}^\mu \lp{\ln\fa t{a_{j+1}T_1}}^{bp+1}
\qquad \mbox{for }t\ge a_{j+1}T_1. 
\end{align*} 
Based on the fact above, we define the sequences $b_j$ and $C_j$ for $j = 0, 1, 2, \cdots$ by
\begin{eqnarray}
b_{j+1} = pb_j+1, \quad & C_{j+1} = \dsp\frac{A_1C_R\lp{\fa 23}^\mu C_j^p}{2(pb_j+1)2^{j+1}}  \label{e8}\\
b_0 = b, \quad & C_0 = A_0.  \label{e9}
\end{eqnarray}
Solving \eqref{e8} and \eqref{e9}, we obtain 
\[
b_j = p^j\lp{b+\fa{1}{p-1}}-\fa{1}{p-1},
\]
and thus,
\[
C_{j+1} = \frac{A_1C_R\lp{\fa 23}^\mu C_j^p}{2b_{j+1}2^{j+1}} \ge \lp{b+\fa{1}{p-1}}^{-1}\frac{A_1C_R\lp{\fa 23}^\mu C_j^p}{2p^{j+1}2^{j+1}}.
\]
Then, 
\begin{align*}
C_j &\ge \frac{BC_{j-1}^p}{2^jp^j}
\ge \cdots\cdots \ge \fa{B^{1+p+p^2+\cdots+p^{j-1}}}{(2p)^{j+p(j-1)+p^2(j-2)+\cdots+p^{j-1}}}C_0^{p^j}, \\
\intertext{and }
\ln C_j &\ge 
\frac{\ln B}{p-1}(p^j-1)-p^j \sum_{k=0}^j 
\frac{k}{p^k}\ln(2p)+p^j\ln C_0, 
\end{align*}
where $B=\lp{b+\fa{1}{p-1}}^{-1}\frac{A_1C_R\lp{\fa 23}^\mu}{2}$. 
For sufficiently large $j$, we have
\[
C_j \ge \exp(Ep^j),
\] 
where
\[
E = \frac{1}{p-1}\min\left(0, \ln B \right) - \sum_{k=0}^{\infty}\frac{k}{p^k}\ln(2p)+\ln A_0. 
\]
Thus, since $F(t)\ge C_j\lp{\ln\fa t{a_jT_1}}^{b_j}
$ holds for $t\ge a_jT_1$, we obtain
\[
F(t) \ge \lp{\ln\fa t{a_jT_1}}^{-1/(p-1)} 
 \exp\left[\left\{E+\lp{b+\fa{1}{p-1}}\ln\lp{\ln\fa t{a_jT_1}}\right\}p^j\right]
  \qquad \mbox{for }t\ge a_jT_1. 
\]
Note that $1\le a_j\le 2$ for all $j=0,1,2,\cdots$. 
We hence have
\begin{equation}
 F(t) \ge \lp{\ln\fa t{T_1}}^{-1/(p-1)} 
 \exp\left[\left\{E+\lp{b+\fa{1}{p-1}}\ln\lp{\ln\fa t{2T_1}}\right\}p^j\right]
  \qquad \mbox{for }t\ge 2T_1. 
 \label{10}
\end{equation}
By choosing $t$ large enough, we can find a positive $\delta$ such that
\[
E+\lp{b+\fa{1}{p-1}}\ln\lp{\ln\fa t{2T_1}}\ge \delta > 0. 
\]
It then follows from \eqref{10}  that $F(t) \lr \infty$ as $j \to \infty$ for 
sufficiently large $t$. We therefore see that the lifespan $T$ of $F(t)$ has to satisfy
\[
T<\exp\lp{CA_0^{-(p-1)/\{b(p-1)+1\}}},  
\]
where $C$ is a constant depending on $R,A_1,\mu,p$ and $b$. 
This completes the proof of the lemma. 
\hfill\qed

\vspace{1cm}

We now prove Theorem 1.1 in the critical case by applying Lemma 3.1. 
Set 
\[
F(t)=\int u(t,x)dx.
\] 
From \eqref{Ho} and \eqref{F2}, 
we have 
\begin{equation}
\int |u|^pdx \ge \fa{C\ep^p}{t^{2}} \qquad \mbox{for }t\ge 1
\label{up}
\end{equation}
since $p=p_F(n(1-\al))$. 
Using \eqref{Fup} and \eqref{up} implies 
\begin{align}
F(t)&\ge C\ep^p \int_1^t \tau^{-\mu}\int_1^\tau s^{\mu-2}dsd\tau \nonumber \\
&\ge C\ep^p \int_1^t \tau^{-\mu-2}
  \int_1^\tau (s-1)^\mu dsd\tau \nonumber\\
&\ge C\ep^p  \int_2^t \tau^{-1}d\tau,  \nonumber
\intertext{therefore, we obtain}
F(t)&\ge C\ep^p \ln\fa t2
\qquad \mbox{for }t\ge 2. 
\label{an}
\end{align}
Finally, by \eqref{Ho} and \eqref{an}, applying Lemma 3.1 with $b=1$ and $A_0=C\ep^p$ 
, we obtain the desired results. 
This completes the proof of Theorem 1.1. 
\hfill\qed \\

\vspace{1cm}

We summarize the upper bounds of the lifespan of blow-up solutions in the subcritical cases including \cite{TW1}.  
We have proved in the present paper that if $1<p<p_F(n(1-\al))$, then 
the upper bound of the lifespan is heatlike:
\begin{equation}
T_{\varepsilon} \le 
C \varepsilon^{- \frac{p-1}{2-n(1-\alpha)(p-1)}}. 
\label{upb3}
\end{equation}
The following two results are shown in \cite{TW1}: \\
If $1<p<p_c(n,\al,\mu)$, then 
the upper bound of the lifespan is wavelike:
\begin{equation}
T_\ep\le C\ep^{\fa{-2p(p-1)}{(1-\al)\gamma(n,p,\al,\mu)}}.
\label{upb1}
\end{equation}  
If $p<1+2/(n(1-\al)+\mu-1)$, then
\begin{equation}
T_\ep\le C\ep^{-\fa{p-1}{2-\{n(1-\al)+\mu-1\}(p-1)}}.  
\label{upb2}
\end{equation}
We recall that $p_c(n,\al,\mu)$ is the positive root of the equation $\gamma(n,p,\al,\mu)=0$ given in \eqref{gm}. 
We note that $2p(p-1)/\{(1-\al)\gamma(n,p,\al,\mu)\}=(p-1)/[2-\{n(1-\al)+\mu-1\}(p-1)]$ yields $p= 2(1-\al)/(n(1-\al)+\mu-1)$, and also that 
$2p(p-1)/\{(1-\al)\gamma(n,p,\al,\mu)\}=(p-1)/\{2-n(1-\al)(p-1)\}$ yields $p= 2(1-\al)/\{n(1-\al)-\mu+1\}$.

Figure 1 below shows the regions of blow-up conditions in the case $n=2$ and $\al=0.6$.  
Note that if $p_c(n,\al,\mu)=p_F(n(1-\al))$, then 
\[
\mu=\mu^\ast\equiv \fa{(1-\al)^2n^2+(1-\al)(1+2\al)n+2}{n(1-\al)+2}. 
\]
We easily see that $\mu^\ast >1$. 

Among the three upper bounds above, if $1<p\le 2(1-\al)/(n(1-\al)+\mu-1)$, then \eqref{upb2} is the best. This is Region (A) shown in Figure 1. 

If $\max\{2(1-\al)/(n(1-\al)+\mu-1), \; 2(1-\al)/(n(1-\al)-\mu+1), \; 1\}<p<p_c(n,\al,\mu)$, then \eqref{upb1} is the best (Region (B) in Figure 1). 

On the other hand, 
if $1<p\le 2(1-\al)/(n(1-\al)-\mu+1)$ and $p<p_F(n(1-\al))$, then 
\eqref{upb3} is the best (Region (C) in Figure 1). 

\begin{figure}
\includegraphics[width=13cm, bb=00 400 500 760, clip]{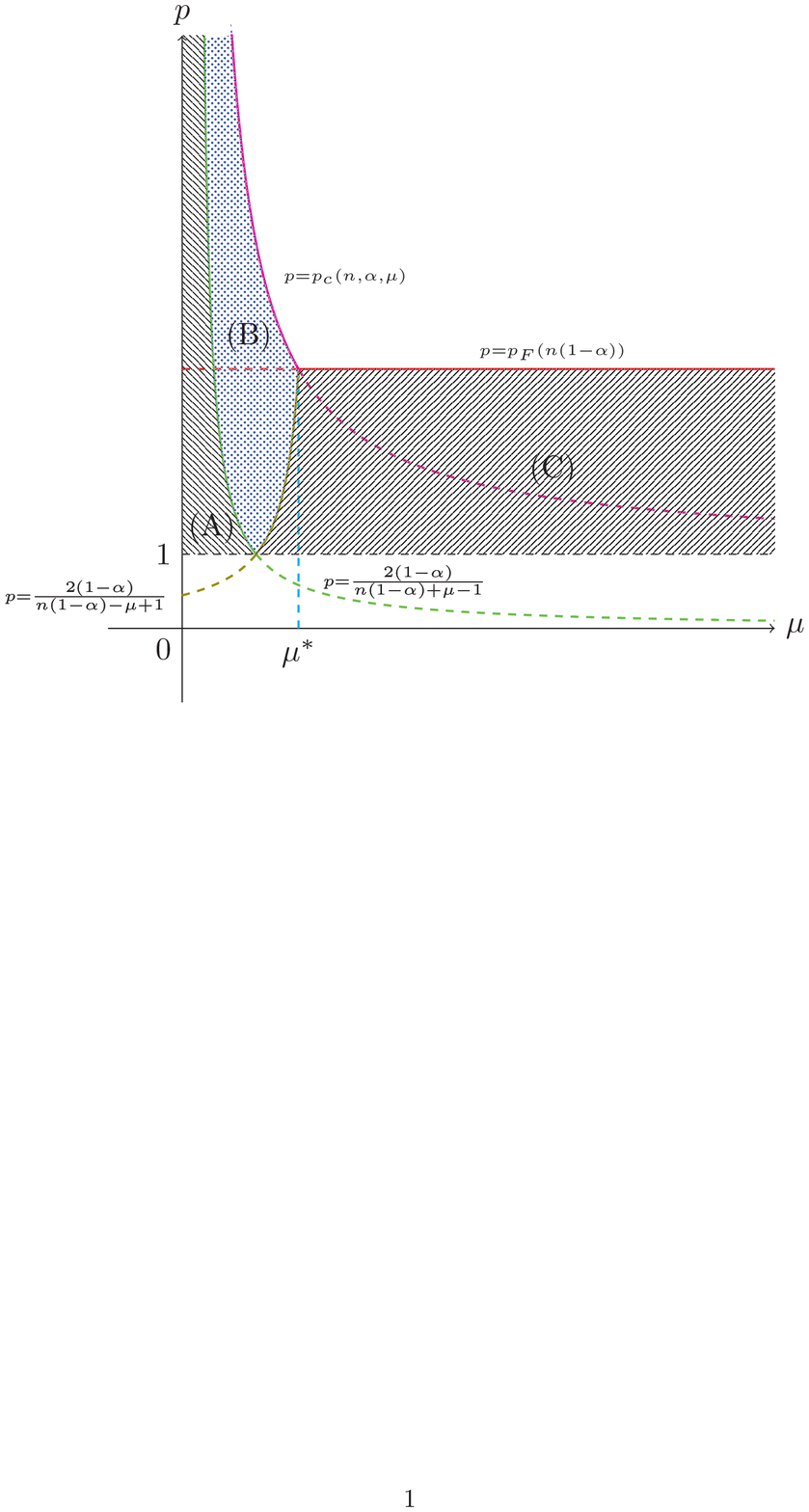}
\caption{Range of blow-up conditions in case $n=2$ and $\alpha=0.6$}

\includegraphics[width=15cm, bb=00 500 500 800, clip]{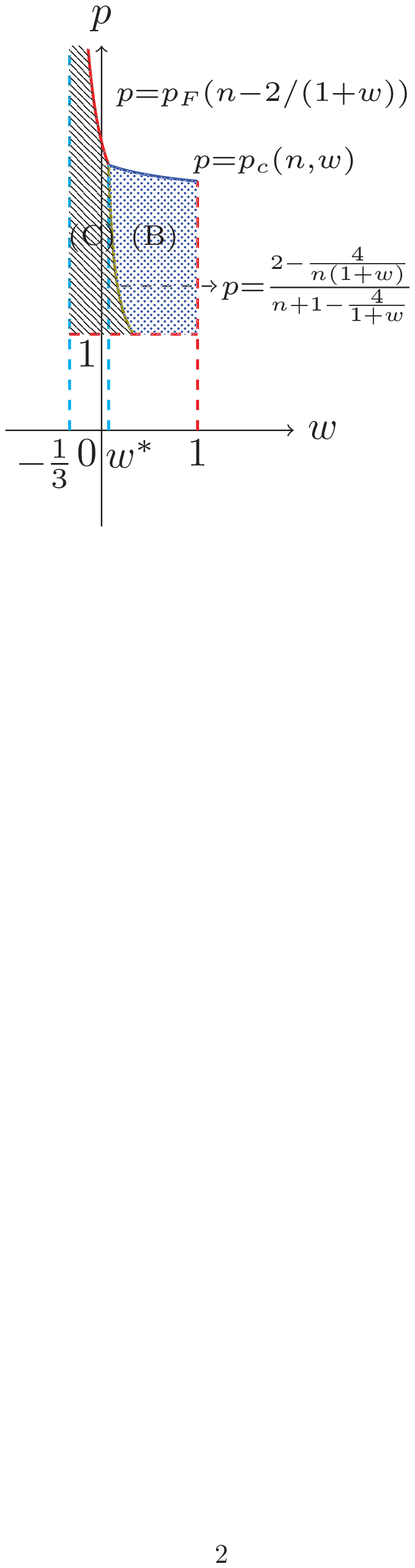}
\caption{Range of blow-up conditions in case $n=3$}

\end{figure}

\section{Wave eq in FLRW. }

\setcounter{equation}{0}

\quad  
We now apply Theorem 1.1 to the original equation \eqref{ore}, which corresponds to \eqref{Prob0} with $\al=2/(n(1+w))$ and $\mu=2/(1+w)$. 
We treat the case $2/n-1<w\le 1$ and $n\ge 2$ so that $1/n\le \al< 1$ and $\mu\ge 1$ in \eqref{Prob0}. 
Observe that this case corresponds to a decelerating expanding universe. 

Assume that $n\ge 2, \; 2/n-1<w\le 1$. 
Then applying Theorem 1.1 to \eqref{ore}, we obtain the following upper bounds of the lifespan: 
\begin{align}
&T_\ep\le C\ep^{\fa{-(p-1)}{2-(n-2/(1+w))(p-1)}}  & &  \mbox{if } 1<p<p_F(n-2/(1+w)), 
\label{31} \\
&T_\ep\le \exp(C\ep^{-(p-1)}) & & \mbox{if } p=p_F(n-2/(1+w))=1+2/\{n-2/(1+w)\}. 
\nonumber
\end{align}
From this result, we see that the blow-up range of $p$ in the flat FLRW spacetime with a decelerating scale factor is larger than that in the Minkowski spacetime because $p_F(n)<p_F(n-2/(1+w))$. 
Moreover, in the subcritical case $p<p_F(n-2/(1+w))$, the lifespan of the blow-up solutions in the FLRW spacetime is shorter than that in the Minkowski spacetime since
$\ep^{- \frac{p-1}{2-(n-2/(1+w))(p-1)}}
<\ep^{-\frac{p-1}{2-n(p-1)}}$ for sufficiently small $\ep$. Therefore, we can say in the heatlike and subcritical case that finite time blow-up  can occur more easily in the FLRW spacetime. The same holds true for the wavelike and subcritical case. 

Indeed, as in \cite{TW1}, substituting $\al=2/(n(1+w))$ and $\mu=2/(1+w)$ into \eqref{upb3} and \eqref{upb1}, we obtain the following results: 
\begin{align}
&T_\ep\le C\ep^{\fa{-2p(p-1)}{\gamma_0(n,p,w)}}  &&  \mbox{if } 1<p<p_c(n,w),  
\label{32}\\
&T_\ep\le \exp(C\ep^{-p(p-1)}) &&  \mbox{if } p=p_c(n,w)>p_F(n-2/(1+w))=1+2/\{n-2/(1+w)\}, \nonumber
\end{align}
where $p_c(n,w)$ is the positive root of the equation $\gamma_0(n,p,w)=0$, and 
\begin{align*}
\gamma_0(n,p,w)&=\lp{1-\fa 2{n(1+w)}}\gamma\lp{n,p,\fa 2{n(1+w)},\fa 2{1+w}}\\
&=-(n-1)p^2+\lp{n+1+\fa 4{n(1+w)}}p+2-\fa 4{n(1+w)}. 
\end{align*}
Recall that $\gamma(n,p,\al,\mu)$ is given in \eqref{gm}.  
Note that $p_c(n,w)>p_S(n)$ and $\gamma_0(n,p,w)>\gamma_S(n,p)
$, where $p_S(n)$ is the Strauss exponent and $\gamma_S(n,p)$ is given in \eqref{snp}.

Figure 2 shows the range of blow-up conditions in terms of $w$ and $p$ in the case $n=3$. 
Note that if $p_F(n-2/(1+w))=p_c(n,w)$, then $w$ is the larger root $w^\ast$ of the equation
\[
n(n^2+n+2)w^2+2n(n-1)^2w+n^3-5n^2+8n-8=0. 
\]
In Region (B) \eqref{32} is better than \eqref{31}, on the other hand, this relation becomes reverse in Region (C). 
We note that Region (A) in Figure 1 does not appear in Figure 2 since $\mu\ge 1$ in our case yields $p> 2(1-\al)/(n(1-\al)+\mu-1)$. 

The remaining case $-1<w\le 2/n-1$ is treated in \cite{TW2}.


\end{document}